\documentclass[openacc]{rstransa}

\usepackage{psfrag}

\newcommand{\RR}{\mathbb{R}}

\newcommand{\ZZ}{\mathbb{Z}}
\newcommand{\PP}{\mathbb{P}}
\newcommand{\EE}{\mathbb{E}}

\newcommand{\Var}{\mathbb{V}{\rm ar}}

\newcommand{\dive}{\operatorname{div}}

\newcommand{\dis}{\displaystyle}

\begin{document}

\title{Some variance reduction methods for numerical stochastic homogenization}
\author{X. Blanc$^1$, C. Le Bris$^2$ and F. Legoll$^2$}
\address{$^1$ Univ. Paris Diderot, Sorbonne Paris Cit\'e, Laboratoire Jacques-Louis Lions, UMR 7598, UPMC, CNRS, F-75205 Paris, France \\$^2$Ecole des Ponts and INRIA, 6 \& 8 avenue Blaise Pascal, 77455 Marne-la-Vall\'ee Cedex 2, France.}

\subject{35R60, 35B27, 65C05}



\corres{Le Bris C.\\
\email{lebris@cermics.enpc.fr}}

\begin{abstract} 
We overview a series of recent works devoted to variance reduction techniques for numerical stochastic homogenization. Numerical homogenization requires solving a set of problems at the micro scale, the so-called corrector problems. In a random environment, these problems are stochastic and therefore need to be repeatedly solved, for several configurations of the medium considered. An empirical average over all configurations is then performed using the Monte-Carlo approach, so as to approximate the effective coefficients necessary to determine the macroscopic behavior. Variance severely affects the accuracy and the  cost of such computations. Variance reduction approaches, borrowed from other contexts of the engineering sciences, can be useful. Some of these variance reduction techniques are presented, studied and tested here. 
\end{abstract}


%
%
%
%
%
%
%
%
%
%
%


\maketitle

\section*{Introduction}

We overview a series of recent works related to some {\em random} multiscale problems motivated by practical problems in Mechanics. For simplicity, we argue on a linear elliptic scalar equation in divergence form:
\begin{equation}
\label{eq:pb0}
\left\{
\begin{array}{l}
\dis -\hbox{\rm div}\left[ A\left(\frac{x}{\varepsilon}\right) \nabla u^\varepsilon \right] = f \quad \text{in} \quad {\cal D}, 
\\ \noalign{\vskip 3pt}
u^\varepsilon=0 \quad \text{on} \quad \partial {\cal D},
\end{array}
\right.
\end{equation}
although the scope of the techniques we describe go beyond this simple setting. The matrix coefficient~$A$ is assumed random stationary. The purpose is to compute the homogenized matrix coefficient~$A^\star$ present in the homogenized equation
\begin{equation}
\label{eq:pb0-star}
\left\{
\begin{array}{l}
-\hbox{\rm div}\left[ A^\star \nabla u^\star \right] = f \quad \text{in} \quad {\cal D}, 
\\ \noalign{\vskip 3pt}
u^\star=0 \quad \text{on} \quad \partial {\cal D},
\end{array}
\right.
\end{equation}
which captures the average behavior of the solution~$u^\varepsilon$ to~\eqref{eq:pb0}.

We begin by recalling in Section~\ref{sec:recall} the basics of homogenization theory, both in the deterministic (periodic) context and in the random context, which are useful for our exposition. Next, we successively present three different variance reduction techniques for the problem considered. We begin in Section~\ref{sec:antithetic} with the classical, general purpose technique of \emph{antithetic variables}. The efficiency of that technique is substantial, but is also limited in particular because the technique does not exploit much the specifics of the problem considered. We present in Section~\ref{sec:controlvariate} the technique of \emph{control variate}, which requires a better knowledge of the problem at hand. A problem simpler to simulate and close to the original problem, in a sense that is made precise below, has to be considered and concurrently solved. The technique uses that knowledge to, effectively, get a much better reduction of the variance. In Section~\ref{sec:sqs}, we expose a slightly different approach, imported from solid state physics, namely that of \emph{special quasi-random structures}. It consists in selecting (somewhat in the spirit of another well-known technique, {\em stratified sampling}) some configurations of the random environment that are more suitable than generic configurations to compute the empirical averages, so as to again minimize the variance. Our final Section~\ref{sec:future} presents some further research directions. 

\medskip

Before we proceed, we mention that we will assume throughout our text that the reader is reasonably familiar with the homogenization theory. We refer to the classical textbooks~\cite{blp,jikov} for this topic. We also mention~\cite{singapour,enumath,smai2013} for general presentations and overviews of the issues examined here, along with some related issues.

\section{Brief overview of classical homogenization settings}
\label{sec:recall}

\subsection{Periodic homogenization}
\label{ssec:per}

To begin with, we recall some well known, basic ingredients of elliptic homogenization theory in the periodic setting. We consider~\eqref{eq:pb0} in a regular, bounded domain ${\cal D}$ in $\RR^d$, and choose the matrix coefficient $A=A_{per}$ to be symmetric and $\ZZ^d$-periodic. Note that, throughout our text, we manipulate for simplicity \emph{symmetric} matrices, but our discussion in Sections~\ref{sec:antithetic} through~\ref{sec:sqs} carries over to non symmetric matrices up to slight modifications.

The corrector problem associated to~\eqref{eq:pb0} in the periodic case reads, for $p$ fixed in $\RR^d$,
\begin{equation}
\label{eq:cor-per-intro}
\left\{
\begin{array}{l}
-\hbox{\rm div}\left(A_{per}(y)\left(p+ \nabla w_p \right)\right) =0, 
\\ \noalign{\vskip 3pt}
\text{$w_p$ is $\ZZ^d$-periodic}.
\end{array}
\right.
\end{equation}
It has a unique solution up to the addition of a constant. Then, the homogenized coefficient in~\eqref{eq:pb0-star} reads
$$
[ A^\star ]_{ij} = \int_Q e_i^T A_{per}(y) \left( e_j + \nabla w_{e_j}(y) \right) dy, 
$$
where $Q=(0,1)^d$ is the unit cube and $(e_i)_{1 \leq i \leq d}$ is the canonical basis of $\RR^d$. Equivalently, $A^\star$ satisfies
$$
\forall p \in \RR^d, \quad A^\star \, p = \int_Q A_{per}(y) \left( p + \nabla w_p(y) \right) dy.
$$
The main result of periodic homogenization theory is that, as $\varepsilon$ goes to zero, the solution $u^\varepsilon$ to~\eqref{eq:pb0} converges to $u^\star$ solution to~\eqref{eq:pb0-star}. The convergence holds in $L^2({\cal D})$ and weakly in $H^1_0({\cal D})$. The correctors $w_{e_i}$ may also be used to ``correct'' $u^\star$ in order to identify  the behavior of $u^\varepsilon$ in the strong topology of $H^1_0({\cal D})$. 

\medskip

Practically, at the price of only computing $d$ periodic problems~\eqref{eq:cor-per-intro}, the solution to problem~\eqref{eq:pb0} can therefore be efficiently approached for $\varepsilon$ small. 

\subsection{Stochastic homogenization}
\label{ssec:random}

Because this is well known and for the sake of brevity, we skip all technicalities related to the definition of the probabilistic setting (we refer e.g. to~\cite{singapour} for all details). We assume that~$A$ is stationary in the sense 
\begin{equation}
\label{eq:stationnarite-disc}
\forall k \in \ZZ^d, \quad A(x+k, \omega) = A(x,\tau_k\omega) \quad \mbox{almost everywhere in $x$, almost surely}
\end{equation}
(where~$\tau$ is an ergodic group action). We consider the boundary value problem~\eqref{eq:pb0} for $A=A(\cdot,\omega)$. Standard results of stochastic homogenization~\cite{blp,jikov} apply and allow to find the homogenized problem. These results generalize the periodic results recalled in Section~\ref{ssec:per}. The solution $u^\varepsilon(\cdot,\omega)$ to~\eqref{eq:pb0} converges to the solution to~\eqref{eq:pb0-star} where the homogenized matrix is now defined as
$$
[A^\star]_{ij} = \EE\left(\int_Q e_i^T A\left(y,\cdot\right)\,\left( e_j+\nabla w_{e_j}(y,\cdot)\right) \,dy\right),
$$
where, for any $p\in \RR^d$, $w_p$ is the solution (unique up to the addition of a random constant) to
\begin{equation}
\label{eq:correcteur-random}
\left\{
\begin{array}{l}
-\hbox{\rm div}\left[A\left(y,\omega\right)\left(p+ \nabla w_p(y,\omega) \right)\right] =0 \quad \mbox{a.s. on $\RR^d$},
\\ \noalign{\vskip 3pt}
\nabla w_p \quad \mbox{is stationary in the sense of~\eqref{eq:stationnarite-disc}}, 
\\ \noalign{\vskip 3pt}
\displaystyle
\EE\left(\int_Q \nabla w_p(y,\cdot)\,dy\right) = 0.
\end{array}
\right.
\end{equation}
Note that $u^\varepsilon$ is a random function, while its homogenized limit $u^\star$ is deterministic since $A^\star$ is deterministic.

A striking difference between the stochastic setting and the periodic setting can be observed comparing~\eqref{eq:cor-per-intro} and~\eqref{eq:correcteur-random}. In the periodic case, the corrector problem is posed on a bounded domain (namely, the periodic cell~$Q$), since the corrector $w_p$ is periodic. In sharp contrast, the corrector problem~\eqref{eq:correcteur-random} of the random case is posed on the whole space $\RR^d$, and cannot be reduced to a problem posed on a bounded domain. The fact that the random corrector problem is posed on the entire space has far reaching consequences for numerical practice. Truncations of problem~\eqref{eq:correcteur-random} have to be considered. The actual homogenized coefficients are only captured in the asymptotic regime. 

More precisely, the deterministic matrix $A^\star$ is usually approximated by the random matrix $A^\star_N(\omega)$ defined by
\begin{equation}
\label{eq:AstarN}
\forall p \in \RR^d, \quad A^\star_N(\omega) \ p = \frac{1}{|Q_N|} \int_{Q_N} A(\cdot,\omega) \left( p + \nabla w^N_p(\cdot,\omega) \right),
\end{equation}
which is obtained by solving the corrector problem on a \emph{truncated} domain, say the cube $Q_N = (0,N)^d$:
\begin{equation} 
\label{eq:correcteur-random-N}
-\dive \Big[ A(\cdot,\omega) \left( p +  \nabla w_p^N(\cdot,\omega)\right) \Big] = 0,
\qquad
w_p^N(\cdot,\omega) \ \mbox{is $Q_N$-periodic}.
\end{equation}
Although $A^\star$ itself is a deterministic quantity, its practical approximation $A^\star_N$ is random. It is only in the limit of infinitely large domains $Q_N$ that the deterministic value is attained. As shown in~\cite{bourgeat2}, we have 
$$
\lim_{N \to \infty} A^\star_N(\omega) = A^\star \quad \text{almost surely.} 
$$
As usual in the random context, the error $A^\star - A^\star_N(\omega)$ may be expanded as
\begin{equation}
\label{eq:error-decomposition}
A^\star - A_N^\star(\omega)
=
\Big(
A^\star - \EE \left[ A_N^\star \right]
\Big)
+
\Big(
\EE \left[ A_N^\star \right]
- A_N^\star(\omega)
\Big),
\end{equation}
that is the sum of a {\em systematic} error and of a {\em statistical} error (the first and second terms in the above right-hand side, respectively). 

A standard technique to compute an approximation of $\EE \left[A^\star_N \right]$ is to consider $M$ independent and identically distributed realizations of the field $A$, solve for each of them the corrector problem~\eqref{eq:correcteur-random-N} (thereby obtaining from~\eqref{eq:AstarN} i.i.d. realizations $A^{\star,m}_N(\omega)$, $1 \leq m \leq M$) and compute the Monte Carlo approximation
$$
\EE \left[ \left( A^\star_N \right)_{ij}\right]
\approx
I^{\rm MC}_M(\omega) := \frac{1}{M} \sum_{m=1}^M \left( A^{\star,m}_N(\omega) \right)_{ij}.
$$
In view of the Central Limit Theorem, we know that $\EE \left[ \left( A^\star_N \right)_{ij} \right]$ asymptotically lies within the confidence interval
$$
\left[
I^{\rm MC}_M - 1.96 \frac{\sqrt{\Var \left[ \left( A^\star_N \right)_{ij} \right]}}{\sqrt{M}}
,
I^{\rm MC}_M + 1.96 \frac{\sqrt{\Var \left[ \left( A^\star_N \right)_{ij} \right]}}{\sqrt{M}}
\right]
$$
with a probability equal to 95 \%. 

For simplicity, and because this is overwhelmingly the case in the numerical practice, we have considered in~\eqref{eq:correcteur-random-N}  \textit{periodic} boundary conditions. These will be the conditions we adopt throughout our study. Other boundary conditions, or approximations, may be employed. The specific choice of approximation technique is motivated by considerations about the decrease of the systematic error in~\eqref{eq:error-decomposition}. Several recent mathematical studies by A.~Gloria and F.~Otto~\cite{gloria-otto} have clarified this issue. The variance reduction techniques we present in this article can be applied to all types of boundary conditions.

\section{Variance reduction using antithetic variables}
\label{sec:antithetic}

We present here a first attempt~\cite{cedya,mprf,banff} to reduce the variance in stochastic homogenization. The technique used for variance reduction is that of \emph{antithetic variables}. 

\medskip

The variance reduction technique using antithetic variables consists in concurrently considering two sets of configurations for the random material instead of only one set. The two sets of configurations will be deduced one from the other. Indeed, fix $M=2 \mathcal{M}$. Suppose that we give ourselves $\mathcal{M}$ i.i.d. copies $\left(A^m(x,\omega)\right)_{1 \leq m \leq \mathcal{M}}$ of $A(x,\omega)$. Construct next $\mathcal{M}$ i.i.d. \emph{antithetic} random fields 
$$
B^m(x,\omega)
= 
T \left( A^m(x,\omega) \right),
\quad
1 \leq m \leq \mathcal{M},
$$
from the $\left(A^m(x,\omega)\right)_{1 \leq m \leq \mathcal{M}}$. The map $T$ transforms the random field $A^m$ into another, so-called \emph{antithetic}, field $B^m$. The transformation is performed in such a way that, for each $m$, $B^m$ has the same law as $A^m$, namely the law of the matrix~$A$. Somewhat vaguely stated, if $A$ was obtained in a coin tossing game (using a fair coin), then $B^m$ would be {\it head} each time $A^m$ is {\it tail} and vice versa. Then, for each $1 \leq m \leq \mathcal{M}$, we solve two corrector problems. One is associated to the original field $A^m$, the other one is associated to the antithetic field $B^m$. Using its solution $v_p^{N,m}$, we define the \emph{antithetic homogenized matrix} $B^{\star,m}_N$, the elements of which read, for any $1 \leq i,j \leq d$, 
$$
\left[ B^{\star,m}_N(\omega) \right]_{ij} = 
\frac{1}{|Q_N|} \int_{Q_N} e_i^T B^m(\cdot,\omega) \ \left(e_j + \nabla v_{e_j}^{N,m}(\cdot,\omega) \right).
$$
And we finally set, for any $1 \leq m \leq \mathcal{M}$,
$$
\widetilde{A}^{\star,m}_N(\omega) := \frac12 \left( A^{\star,m}_N(\omega) +  B^{\star,m}_N(\omega) \right).
$$
Since $A^m$ and $B^m$ are identically distributed, so are $A^{\star,m}_N$ and $B^{\star,m}_N$. Thus, $\widetilde{A}^{\star,m}_N$ is unbiased (that is, $\mathbb{E}\left(\widetilde{A}^{\star,m}_N\right)=\mathbb{E}\left({A}^{\star,m}_N\right)$). In addition, it satisfies:
$$
\widetilde{A}^{\star,m}_N \underset{N \rightarrow + \infty}{\longrightarrow} A^\star \ \mbox{almost surely},
$$
because $A^m$ and $B^m$ are ergodic. The hope is that the new approximation~$\widetilde{A}^{\star,m}_N$ has less variance than the original one~$A^{\star,m}_N$. It is indeed the case under appropriate assumptions. 

\medskip

The approach has been studied theoretically in~\cite{cedya,mprf,banff}, in the  one-dimensional setting and in some specific higher dimensional cases. The approach is shown to qualitatively reduce the variance. A quantitative assessment of the reduction is however out of reach. Only numerical tests can provide some information in this direction.

\medskip

The tests we have performed in~\cite{cedya,banff} concern various ``input'' random fields $A(\cdot,\omega)$, some i.i.d.,  some correlated, with various correlation lengths. In these settings, we have investigated variance reduction on a typical diagonal~$\left[A^\star_N(\omega) \right]_{11}$, or off-diagonal~$\left[A^\star_N(\omega) \right]_{12}$ entry of the approximate homogenized matrix~$A^\star_N(\omega)$, as well as on the eigenvalues of the matrix $A^\star_N(\omega)$, and the eigenvalues of the associated differential operator~$L = - \hbox{\rm div} \left[ A^\star_N(\omega) \nabla \cdot \right]$ (supplied with homogeneous Dirichlet boundary conditions on~$\partial{\cal D}$). 

\begin{figure}[htbp]
\centering\includegraphics[width=3.5in]{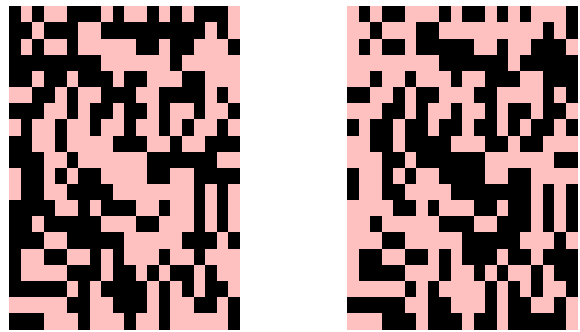}
\caption{An example of realization of $A(x,\omega)$ together with its antithetic field $B(x,\omega)$ (reproduced from~\cite{cedya}).\label{fig_B12}}
\end{figure}

Let us give one such example. Consider, in dimension two, the matrix $A(x,\omega)$ defined by
\begin{equation}
\label{eq:1}
A(x,\omega) = \sum_{k\in\ZZ^2} {\bf 1}_{Q+k}(x) \, a_k(\omega)
\begin{pmatrix}
1 & 0 \\ 0 & 1
\end{pmatrix}
,
\end{equation}
where $Q=(0,1)^2$, $\left(a_k\right)_{k\in\ZZ^2}$ is an i.i.d sequence of random Bernoulli variables such that $\PP(a_k = \alpha) = \PP(a_k = \beta) = 1/2$, with $\alpha = 3$ and $\beta = 20$. An example of the realization of each matrix field $A(x,\omega)$ and $B(x,\omega)$ is given in Figure~\ref{fig_B12} (in black, the value $\alpha$ and in pink, the value $\beta$).

\begin{figure}[htbp]
\psfrag{Number of cells: \(2N\)^2}{Size of $Q_N$}
\centering\includegraphics[width=3.5in]{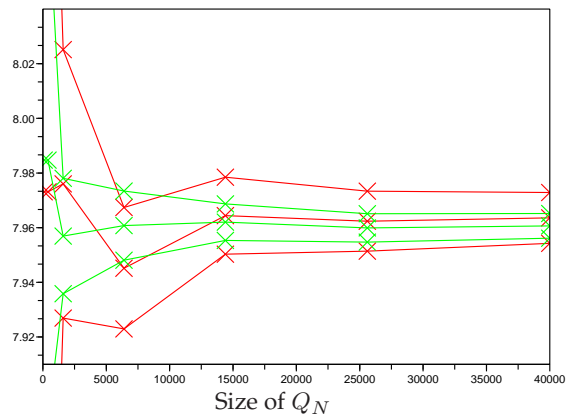}
\caption{Estimation of $A_{11}^\star$ (with confidence interval) with respect to $|Q_N|$ (in red, the classical MC strategy, in green the antithetic variable strategy; reproduced from~\cite{cedya}). \label{fig_B12_result}}
\end{figure}

We then compare two computations with identical cost. For this purpose, we first use a classical Monte Carlo method with $2M$ draws (with here $2M=100$). Second, we apply the antithetic variable technique using only $M$ draws. Since we solve two corrector problems for each of the draws (one for $A^m$ and one for $B^m$), the numerical cost is equal to the cost of the classical computation. The results are shown in Figure~\ref{fig_B12_result}, where we can see that the (numerically estimated) variance is reduced.

\begin{figure}[htbp]
\psfrag{Taille de $Q_N$}{Size of $Q_N$}
\psfrag{CPU time gain}{ \hspace{-5mm} Variance ratio}
\centering\includegraphics[width=2.5in,angle=270]{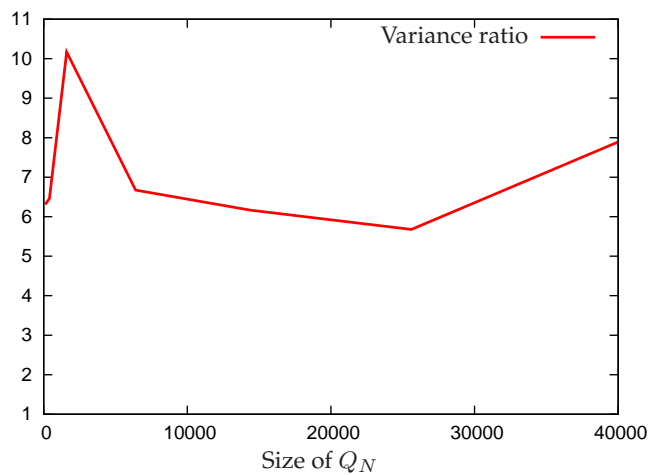}
\caption{Efficiency of the variance reduction (same CPU time, variance ratio).\label{fig_B12_variance}}
\end{figure}

A more precise estimate of the efficiency of the approach is given on Figure~\ref{fig_B12_variance}, in which we have plotted the variance ratio with respect to the size of the computational domain. We see that the gain is not very sensitive to this size, and is at least of about $6$ on this example. This means that, given a computational cost, the approach improves the accuracy by a factor $\sqrt 6 \approx 2.45$. Equivalently, for a given accuracy, the computational cost is reduced by a factor $6$.

\medskip

Our numerical results (see~\cite{cedya,banff} for comprehensive details) show that the technique may be applied to a large variety of situations and has proved efficient whatever the output considered. In addition, we have shown in~\cite{dcds} that this technique carries over to nonlinear stochastic homogenization problems, when the problem at hand is formulated as a variational convex problem. In all the test cases we have considered, variance is systematically reduced. We observed however that the ratio of reduction is not spectacular. This has motivated the consideration of alternative techniques, expected to be (and indeed observed to be) more efficient than the antithetic variables technique.

\section{Control variate technique}
\label{sec:controlvariate}

The \emph{control variate} approach is a variance reduction technique known to be potentially much more efficient than the antithetic variable technique. It however asks to have beforehand a better information on the random quantity of interest that is simulated. In the context of homogenization, the works~\cite{minvielle,minvielle-these} present a first possible investigation of the efficiency of this technique.

The specific setting considered as control variate is a periodic setting slightly perturbed using a random field modeled by a Bernoulli variable which we now briefly describe, before turning to the variance reduction technique itself. 

\subsection{Our specific choice of control variate: a perturbation approach}
\label{ssec:arnaud}

One approach, described in full details in~\cite{cras-arnaud,arnaud1,arnaud2}, addressing the random material as a small perturbation of a periodic material, consists in considering
\begin{equation}
\label{eq:perturb2}
A_\eta(x, \omega) = A_{per}(x) + b_\eta(x, \omega) C_{per}(x),
\end{equation}
where, with evident notation, $A_{per}$ is a $\ZZ^d$-periodic matrix modeling the unperturbed material and $C_{per}$ is a $\ZZ^d$-periodic matrix modeling the perturbation. We take 
$$
b_\eta(x,\omega) = \sum_{k \in \ZZ^d} \mathbf{1}_{Q+ k}(x) B_\eta^k(\omega),
$$
where the $B_\eta^k$ are, say, independent identically distributed scalar random variables. One particularly interesting case (see~\cite{cras-arnaud,arnaud1,arnaud2} for other cases) is when the common law of the $B_\eta^k$ is assumed to be a Bernoulli law of (presumably small) parameter~$\eta$:
$$
\PP(B_\eta^k = 1) = \eta, \qquad \PP(B_\eta^k = 0) = 1-\eta.
$$
A formal approach introduced in the above works (which has subsequently been studied and proved correct in~\cite{mourrat,gloria-arma}) to efficiently perform homogenization in that context starts with observing that, in the corrector problem
%
\begin{equation}
\label{eq:correc-intuition}
-\hbox{\rm div}\left[A_\eta\left(y,\omega\right)\left(p+ \nabla w_p(y,\omega) \right)\right] =0,
\end{equation}
the only source of randomness comes from the coefficient $A_\eta\left(y,\omega\right)$. Therefore, if one knows the law of this coefficient, one knows the law of the corrector function $w_p(y,\omega)$ and therefore may compute the homogenized coefficient~$A^\star_\eta$, the latter being a function of this law. When the law of $A_\eta$ has an expansion in terms of a small coefficient, so has the law of $w_p$. Consequently, $A_\eta^\star$ can be obtained as an expansion. Heuristically, on the cube~$Q_N = [0,N]^d$ and at order 1 in $\eta$, the probability to get the perfect periodic material (entirely modeled by the matrix~$A_{per}$) is $(1-\eta)^{N^d}\approx 1-N^d\eta+O(\eta^2)$, while the probability to obtain the unperturbed material on all cells except one (where $A_\eta = A_{per} + C_{per}$) is $N^d\,(1-\eta)^{N^d-1}\eta\approx N^d\eta+O(\eta^2)$. All other configurations, with two or more cells perturbed, yield contributions of order higher than or equal to $\eta^2$. This gives the intuition (and this intuition can be turned into a mathematical proof) that the first order correction indeed comes from the difference between the material perfectly periodic except on one cell and the perfect material itself: 
\begin{equation}
\label{eq:dev_eta}
A_\eta^\star = A_{per}^\star + \eta A_{1,\star} + o(\eta),
\end{equation} 
where $A_{per}^\star$ is the homogenized matrix for the unperturbed periodic material and
$$
A_{1,\star} = \lim_{N \rightarrow + \infty} A_{1,\star,N}, 
$$
with
\begin{equation}
\label{eq:utile}
A_{1,\star,N} \, e_i = \int_{Q_N}\left[(A_{per}+\mathbf{1}_Q C_{per})(\nabla w_i^N + e_i ) - A_{per}(\nabla w_i^0 + e_i )\right], 
\end{equation}
where $w_i^0$ is the corrector for $A_{per}$ (i.e. the solution to~\eqref{eq:cor-per-intro}), and $w_i^N$ solves
$$
-\mathrm{div}\left((A_{per}+ \mathbf{1}_Q C_{per}) (\nabla w_i^N + e_i ) \right) = 0 \quad  \mathrm{in} \quad Q_N, \quad \text{$w_i^N$ is $Q_N$-periodic}.
$$


The approach has been extensively tested. It is observed that, using the perturbative approach, the large $N$ limit is already very well approached for small values of $N$. The computational efficiency of the approach is clear: solving the two {\em periodic} problems with coefficients $A_{per}$ and $A_{per} + \mathbf{1}_Q C_{per}$ for a limited size $N$ is much less expensive than solving the original, {\em random} corrector problem for a much larger size $N$.

When the second order term is needed, configurations with two defects have to be computed. They all can be seen as a family of PDEs, parameterized by the geometrical location of the defects. Reduced basis techniques have been shown in~\cite{thomines-redbasis} to  allow for a definite speed-up in the computation. 

\subsection{Variance reduction}
\label{ssec:control-variate}

We now again consider the setting defined by~\eqref{eq:perturb2}, except that, now, the parameter~$\eta$ of the Bernoulli law is \emph{not} taken small. The expansion technique employed in Section~\ref{ssec:arnaud} is therefore inaccurate. It can however serve for the construction of a control variate, useful to reduce the variance. 
 
Determining the field $A(x,\omega)$, given by~\eqref{eq:perturb2}, on the truncated domain~$Q_N$ amounts to drawing $B^k_\eta(\omega)$ in each cell $Q+k$ in~$Q_N$. This allows to compute the associated (approximate) homogenized coefficient $A^\star_N(\omega)$ from the solution to the corrector problem~\eqref{eq:correc-intuition} truncated on $Q_N$. In parallel to this task, we \emph{reconstruct} from the specific realization of the set of $B_\eta^k(\omega)$ a field that is used as a control variate. More precisely, we set
\begin{equation}
\label{eq:cv_c}
C^\star_N(\omega) = A^\star_N(\omega) - \rho \left(A_{per}^\star + A_1^{\star,N}(\omega) - \EE \left[ A_{per}^\star + A_1^{\star,N}(\omega) \right] \right).
\end{equation}
In this formula,  
$$
A_1^{\star,N}(\omega)= \frac{1}{|Q_N|} \sum_{k+Q\subset Q_N} B^k_\eta(\omega) \ {\cal A}^{\rm 1 \, def}_k, 
$$
where ${\cal A}^{\rm 1 \, def}_k$ is the deterministic coefficient corresponding to the case of \emph{one} defect located at position~$k$ in $Q_N$ (it is actually independent of~$k$ and equal to $A_{1,\star,N}$ defined by~\eqref{eq:utile}). The parameter~$\rho$ in~\eqref{eq:cv_c} is a deterministic parameter, a classical ingredient of control variate techniques, which is optimized in terms of the estimated variances of the objects at play. It is crucial to note that the expectation of~$A_1^{\star,N}(\omega)$ is analytically computable. Since by construction $\EE\left(C^\star_N\right) = \EE\left(A^\star_N\right)$, the technique then consists in approximating the former (thus the latter) by an empirical mean. The theoretical study and the numerical tests in~\cite{minvielle} show that the variance of $C^\star_N$ is smaller than that of $A^\star_N$, and hence that the quality of the approximation is improved. 

\medskip

\begin{figure}[htbp]
\psfrag{Supercell sidelength N}{\footnotesize \hspace{7mm} $N$}
\psfrag{Homogenized coefficient}{\footnotesize \hspace{-7mm} Homogenized coefficient}
\centering\includegraphics[width=4in]{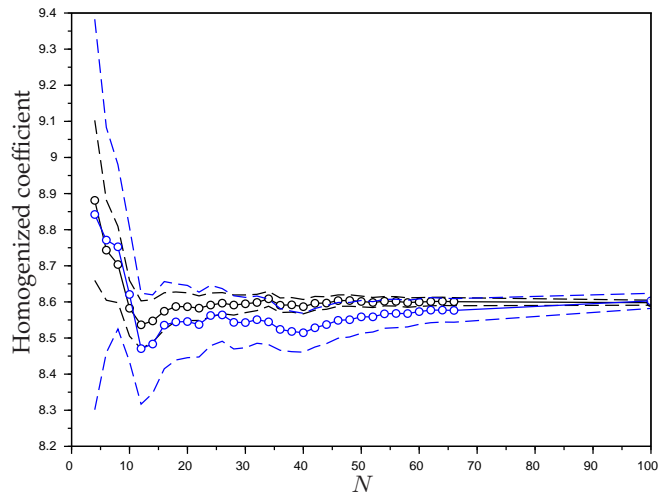}
\caption{Estimation of $A_{11}^\star$ together with its confidence interval (computed using $M=100$ i.i.d. realizations), for the classical MC simulation (in blue) and with the control variate approach~\eqref{eq:cv_c} (in black) (reproduced from~\cite{minvielle}).\label{CV_result}}
\end{figure}

As an illustration, we use a similar case as in Section~\ref{sec:antithetic}, namely~\eqref{eq:1} with $\alpha = 3$ and $\beta = 23$. This case falls within the framework~\eqref{eq:perturb2} with $\eta = 1/2$. This is hence not a perturbative setting. Applying the above strategy based on~\eqref{eq:cv_c} provides the results of Figure~\ref{CV_result}, where the variance is reduced by a factor close to 6, that is, comparable to the technique of antithetic variables.

\begin{figure}[htbp]
\psfrag{Supercell sidelength N}{\footnotesize \hspace{7mm} $N$}
\psfrag{Homogenized coefficient}{\footnotesize \hspace{-7mm} Homogenized coefficient}
\centering\includegraphics[width=4in]{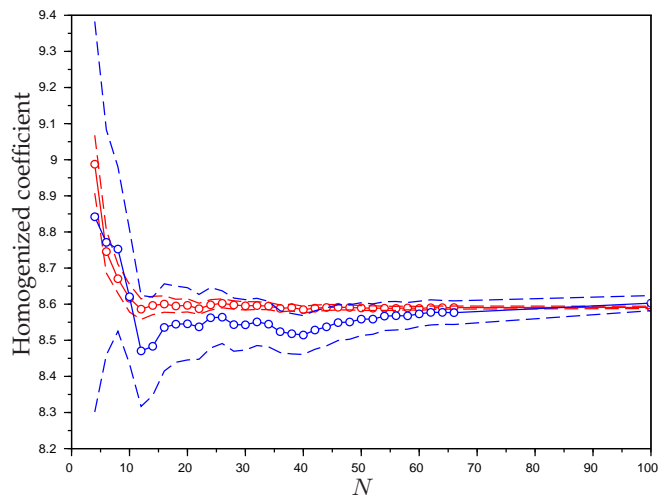}
\caption{Estimation of $A_{11}^\star$ together with its confidence interval (computed using $M=100$ i.i.d. realizations), for the classical MC simulation (in blue) and with the second-order control variate approach (in red) (reproduced from~\cite{minvielle}).\label{CV_result_2}}
\end{figure}

It is also possible to use a second order expansion with respect to $\eta$ in~\eqref{eq:dev_eta}, and include in the control variate both terms, namely the deterministic coefficients corresponding to the case of one and two defects in $Q_N$. Here, one needs additional parameters playing the role of $\rho$ above, in order to ensure substantial variance reduction (see the details in~\cite{minvielle}). The variance reduction of such a case, of the order of 40, is represented on Figure~\ref{CV_result_2}.

\section{Special Quasi-Random Structures}
\label{sec:sqs}

The variance reduction approach we now overview has been originally introduced by other authors for a slightly different purpose in atomistic solid-state science~\cite{vonPezoldDickFriakNeugebauer2010,WeiFerreiraBernardZunger1990,ZungerWeiFerreiraBernard1990}. It carries the name SQS, abbreviation of {\em Special Quasirandom Structures}. The approach has been adapted to the homogenization context in~\cite{sqs,minvielle-these} to which we refer the reader for a more detailed presentation.

\subsection{Motivation and formal derivation of SQS conditions}
\label{sec:derivation_SQS}

In order to convey to the reader the intuition of the original approach, we first consider here a simple one-dimensional setting, which illustrates the difficulties of a generic problem. We consider a linear chain of atomistic sites of two species $A$ and $B$ that interact by a nearest-neighbour interaction potential $V_{AA}$, $V_{AB}$ and $V_{BB}$.

In order to compute the energy per unit particle of that atomistic system, one has to consider all possible such infinite sequences, and for each of them its normalized energy
\begin{equation}
\label{eq:eq_V}
\lim_{N \to \infty} \frac{1}{2N +1} \sum_{i = - N}^N V_{X_{i+1}X_i},
\end{equation}
where $X_i$ denotes the species present at the $i$-th site for that particular configuration. The ``energy'' of the system is then defined as the  \textit{expectation} of~\eqref{eq:eq_V} over all possible configurations. The approach introduced in~\cite{vonPezoldDickFriakNeugebauer2010,WeiFerreiraBernardZunger1990,ZungerWeiFerreiraBernard1990} consists in \textit{selecting} specific truncated configurations $(X_i)_{-N \leq i \leq N}$ of atomic sites that satisfy statistical properties usually obtained only in the limit of infinitely large $N$.

The first such statistical property is the volume fraction, namely the proportion of species $(A, B)$ present on average: one only considers truncated sequences $(X_i)_{-N \leq i \leq N}$ that \textit{exactly} reproduce that volume fraction. Similarly, one may only consider truncated sequences $(X_i)_{-N \leq i \leq N}$ that, in addition to exhibiting the exact volume fraction, have an average energy $\displaystyle \frac{1}{2N+1} \sum_{i = -N}^N V_{X_{i+1}X_i}$ \textit{equal} to $\displaystyle {\cal E} := \frac14 \left(V_{AA} + 2V_{AB} + V_{BB} \right)$. And so on and so forth for other quantities of interest. 

Mathematically, this \textit{selection} of suitable configurations among all possible configurations amounts to replacing the computation of an expectation by that of a \textit{conditional} expectation.

The above simplistic model can of course be replaced by more elaborate models, with more sophisticated quantities to compute, and more demanding statistical quantities to condition the computations with. The bottom line of the approach remains the same, and we now describe its adaptation so as to construct a variance reduction approach for numerical random homogenization.

\medskip

To start with, we assume that the matrix valued random coefficient $A$ present in~\eqref{eq:pb0} reads as
\begin{equation}
\label{eq:A_form}
A_\eta(x, \omega) = C_0(x, \omega) + \eta \ \chi(x, \omega) C_1(x, \omega)
\end{equation}
for some presumably small scalar coefficient $\eta$, and where we assume that $C_0$ and $C_1$ are two stationary, coercive, uniformly bounded matrix fields, that $C_0 - C_1$ is coercive, and that $\chi$ is a stationary scalar field with values in $[-1,1]$. Under these assumptions the matrix $A_\eta$ is stationary, bounded and coercive, uniformly with respect to $\omega$. Since $\eta$ is small, $A_\eta$ is intuitively a perturbation of the matrix valued field $C_0(x, \omega)$.

As already performed above, we may expand all quantities of homogenization theory in powers of the small parameter~$\eta$. In particular, the approximations~$\nabla w^N_\eta$ and $A^{\star,N}_\eta$ of, respectively, the corrector~$\nabla w_\eta$ and the homogenized matrix~$A^\star_\eta$ on the truncated domain $Q_N$, can be expanded in powers of $\eta$:
\begin{eqnarray}
\nonumber
\nabla w^N_\eta(\cdot,\omega) &=& \nabla w^N_0(\cdot,\omega) + \eta \nabla u^N_1(\cdot,\omega) + \eta^2 \nabla u^N_2(\cdot,\omega) + o(\eta^2),
\\
\label{eq:monde_reel}
A^{\star,N}_\eta(\omega) &=& A^{\star,N}_0(\omega) + \eta A^{\star,N}_1(\omega) + \eta^2 A^{\star,N}_2(\omega) + o(\eta^2).
\end{eqnarray}
Inserting these two expansions in~\eqref{eq:correcteur-random-N} and~\eqref{eq:AstarN}, one easily sees that 
$$
\left\{
\begin{array}{ccc}
- \dive C_0 (p + \nabla w^N_0) = 0& \quad \text{in } Q_N, &\quad \text{$w^N_0$ is $Q_N$-periodic},
\\ \noalign{\vskip 3pt}
- \dive C_0 \nabla u^N_1 = \dive \left[ \chi C_1 (p+\nabla w^N_0) \right]& \quad \text{in } Q_N, &\quad \text{$u^N_1$ is $Q_N$-periodic},
\\ \noalign{\vskip 3pt}
- \dive C_0 \nabla u^N_2 = \dive \left[ \chi C_1 \nabla u^N_1 \right]& \quad \text{in } Q_N, &\quad \text{$u^N_2$ is $Q_N$-periodic},
\end{array}
\right.
$$
and that the random variables $A^{\star,N}_0(\omega)$, $A^{\star,N}_1(\omega)$ and $A^{\star,N}_2(\omega)$ read as
\begin{eqnarray*}
A^{\star,N}_0(\omega) \, p &=& \frac{1}{|Q_N|} \int_{Q_N} C_0(\cdot,\omega) (p + \nabla w^N_0(\cdot,\omega)),
\\ \noalign{\vskip 3 pt}
A^{\star,N}_1(\omega) \, p &=& \frac{1}{|Q_N|} \int_{Q_N} \chi(\cdot,\omega) C_1(\cdot,\omega) (p+ \nabla w^N_0(\cdot,\omega)) + \frac{1}{|Q_N|} \int_{Q_N} C_0(\cdot,\omega) \nabla u^N_1(\cdot,\omega),
\\ \noalign{\vskip 3 pt}
A^{\star,N}_2(\omega) \, p &=& \frac{1}{|Q_N|} \int_{Q_N} \chi(\cdot,\omega) C_1(\cdot,\omega) \nabla u^N_1(\cdot,\omega) + \frac{1}{|Q_N|} \int_{Q_N} C_0(\cdot,\omega) \nabla u^N_2(\cdot,\omega).
\end{eqnarray*}

In line with the motivation we have mentioned above in the context of solid state science, we are now in position to introduce the conditions that we use to select particular configurations of the environment within $Q_N$.

For finite fixed $N$, we say that a configuration $\omega$ satisfies the SQS conditions of order up to $k$ if, for any $0 \leq j \leq k$, the coefficient $A^{\star,N}_j(\omega)$ of the expansion~\eqref{eq:monde_reel} exactly matches the corresponding coefficient~$A^\star_j$ of the analogous expansion of the exact homogenized matrix coefficient $A^\star_\eta$. More explicitly, we speak about the SQS condition of 
\begin{itemize}
\item order 0 if $A^{\star,N}_0(\omega) = A^\star_0$, that is to say, for any $p \in \RR^d$,
\begin{equation}
\label{eq:cond0}
\frac{1}{|Q_N|}\int_{Q_N} C_0(x, \omega) (p + \nabla w^N_0(x, \omega)) dx = \EE\left[\int_Q C_0 (p + \nabla w_0)\right],
\end{equation}
\item order 1 if $A^{\star,N}_1(\omega) = A^\star_1$, that is to say, for any $p \in \RR^d$,
\begin{multline}
\label{eq:cond1}
\frac{1}{|Q_N|}\int_{Q_N} \left( \chi(x, \omega)C_1(x, \omega) (p+ \nabla w^N_0(x, \omega)) + C_0(x, \omega) \nabla u^N_1(x, \omega) \right) dx \\ 
= \EE\left[ \int_Q \chi C_1 (p+ \nabla w_0) + C_0 \nabla u_1\right],
\end{multline}
\item order 2 if $A^{\star,N}_2(\omega) = A^\star_2$, that is to say, for any $p \in \RR^d$,
\begin{multline}
\label{eq:cond2}
\frac{1}{|Q_N|}\int_{Q_N} \left( \chi(x, \omega)C_1(x, \omega) \nabla u^N_1(x, \omega) +C_0(x, \omega) \nabla u^N_2 (x, \omega) \right) dx \\ 
= \EE\left[\int_Q\chi C_1 \nabla u_1 + C_0 \nabla u_2\right].
\end{multline}
\end{itemize}

It is easily observed that using such particular configurations that satisfy the SQS conditions of order up to $k$ we have, in the perturbative setting considered here,
\begin{equation}
\label{eq:ordre_k}
A^{\star,N}_\eta(\omega)-A^\star_\eta = o(\eta^k).
\end{equation}
Taking the expectation over such configurations therefore formally provides a more accurate approximation of $A^\star_\eta$. Of course, the purpose is to apply the approach \emph{beyond} the perturbative setting. A property such as~\eqref{eq:ordre_k} cannot be expected any longer since the homogenized matrix $A^\star$ is no longer a series in a small coefficient that encodes a perturbation. Nevertheless, it can be expected that selecting the configurations using these conditions may improve the approximation, in particular by reducing the variance. 

To make the computation of the right-hand sides of the above conditions practical (since in theory they can only be determined using an asymptotic limit, and are therefore as challenging to compute in practice as $A^\star$ itself), we restrict the generality of our setting. We assume that, in~\eqref{eq:A_form}, $C_0(x, \omega) = C_0$ is a deterministic, \textit{constant} matrix, $C_1(x, \omega) = C_1(x)$ is a deterministic, $\ZZ^d$-periodic matrix, and that $\dis \chi(x, \omega) = \sum_{k \in \ZZ^d} X_k(\omega) 1_{k+Q}(x)$, where $X_k(\omega)$ are identically distributed, not necessarily independent, bounded random variables. For the sake of simplicity, we also assume here that 
$$
\EE\left[ X_0 \right] = 0
$$
and refer to~\cite{sqs,minvielle-these} for more general cases. After a tedious but not complicated calculation (the detail of which is provided in~\cite{sqs,minvielle-these}), we obtain that the two conditions~\eqref{eq:cond1}--\eqref{eq:cond2} rewrite as
\begin{eqnarray}
\label{eq:cond12_un}
\dis \frac{1}{|Q_N|} \sum_{k \in \ZZ^d \cap Q_N} X_k(\omega) &=& 0,
\\ \noalign{\vskip 3pt}
\label{eq:cond12_deux}
\dis \frac{1}{|Q_N|} \sum_{k, j \in Q_N \cap \ZZ^d} X_k(\omega) X_j(\omega) I_{k,j}^N &=& \sum_{k \in \ZZ^d} \EE[X_0 X_k] I_k^\infty,
\end{eqnarray}
respectively, where $\dis I_k^\infty = \int_{k+Q} C_1 \nabla \phi_1$ and $\dis I_{k,j}^N = \int_{Q+j} C_1(x) \nabla \phi_1^N(x-k)dx$. In these expressions, $\phi_1$ is the (unique up to the addition of a constant) solution in $\left\{ v \in L^2_{\rm loc}(\RR^d), \ \ \nabla v \in (L^2(\RR^d))^d \right\}$ to
$$
- \dive \left[ C_0 \nabla \phi_1 \right] = \dive \left[ \mathbf{1}_Q C_1p\right] \quad \text{in $\RR^d$},
$$
while $\phi_1^N$ is the (unique up to the addition of a constant) solution to
$$
- \dive \left[ C_0 \nabla \phi_1^N \right] = \dive \left[ \mathbf{1}_{Q}  C_1p \right] \ \ \text{in $Q_N$}, \quad \text{$\phi_1^N$ is $Q_N$-periodic}.
$$
The conditions~\eqref{eq:cond12_un}--\eqref{eq:cond12_deux} are called the SQS 1 and SQS 2 conditions. On the other hand, in the particular setting chosen, condition~\eqref{eq:cond0} (SQS 0, in some sense) is easily seen to be systematically satisfied when $N$ is an integer and the truncated approximation of~\eqref{eq:correcteur-random} that is chosen is the periodic approximation~\eqref{eq:correcteur-random-N}. 

\subsection{Selection Monte Carlo sampling}
\label{sec:algo}

The classical Monte Carlo sampling consists in successively generating a random configuration~$\omega_m$, solving the truncated corrector problem~\eqref{eq:correcteur-random-N} for that configuration, computing~$A^\star_N(\omega_m)$, and finally computing the empirical mean~$\dis {\cal I}^M_{MC} := \frac1M \sum_{m=1}^M A^\star_N(\omega_m)$ as an approximation for~$A^\star$.

In our selection Monte Carlo sampling, we systematically test whether the generated configuration satisfies the required SQS conditions, up to a certain tolerance, and reject it if it does not, \emph{before} solving the corrector problem~\eqref{eq:correcteur-random-N} for that configuration and letting it contribute to the empirical mean.  

In full generality, the cost of Monte Carlo approaches is usually dominated by the cost of draws, and therefore selection algorithms are targeted to reject as few draws as possible. In contrast, in the present context where boundary value problems such as~\eqref{eq:correcteur-random-N} are to be solved repeatedly, the cost of draws for the configuration is negligible compared to the cost of the solution procedure for such boundary value problems. Likewise, evaluating the quantities present in~\eqref{eq:cond12_un}--\eqref{eq:cond12_deux} is inexpensive. Therefore, the purpose of the selection mechanism is to limit the number of boundary value problems to be solved, even though this comes at the (tiny) price of rejecting many configurations. We also note that, as  for any selection procedure, our selection may introduce a bias (i.e. a modification of the systematic error in~\eqref{eq:error-decomposition}). The point is to ensure that the gain in variance dominates the bias introduced by the selection approach. 

We have studied the approach theoretically in~\cite{sqs,minvielle-these}. It is shown therein that the estimator provided (at least the simplest variant of our approach) converges towards the homogenized coefficient $A^\star$ when the truncated domain converges to the whole space. The efficiency of the approach is also theoretically demonstrated for some particular and simple situations (such as the one-dimensional setting). A comprehensive experimental study of the approach has been completed. In particular, since it is often necessary to enforce the desired conditions only up to some tolerance, we have investigated in~\cite{sqs,minvielle-these} how this tolerance affects the quality of the approximation and the efficiency of the approach. We have observed that the approach is robust in this respect.

\begin{figure}[htbp]
\centering\includegraphics[width=3in]{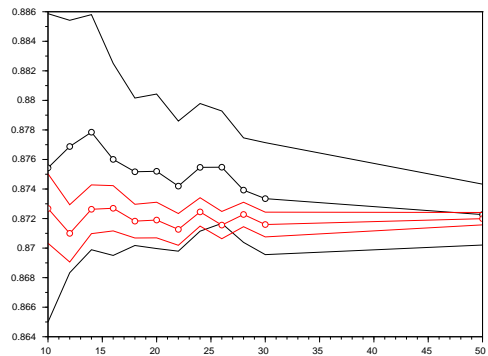}
\caption{Estimation of $A_{11}^\star$ together with its confidence interval (computed using $M=100$ i.i.d. realizations) as a function of $N$, for the classical MC simulation (in black) and with the SQS approach based on~\eqref{eq:cond12_un} (in red) (reproduced from~\cite{sqs}). \label{fig:SQS_un}}
\end{figure}

\begin{figure}[htbp]
\psfrag{Red = SQS 1st order, Blue = SQS 2nd order, Black- = MC, Black-- = exact}{}
\psfrag{meanvalue}{}
\psfrag{supercell size}{\footnotesize \hspace{-2mm} $N$}
\centering\includegraphics[width=3in]{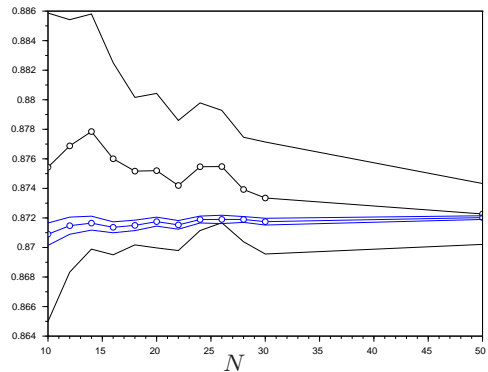}
\caption{Estimation of $A_{11}^\star$ together with its confidence interval (computed using $M=100$ i.i.d. realizations) as a function of $N$, for the classical MC simulation (in black) and with the SQS approach based on~\eqref{eq:cond12_un} and~\eqref{eq:cond12_deux} (in blue) (reproduced from~\cite{sqs}).\label{fig:SQS_deux}}
\end{figure}

We include here a typical illustration of the efficiency of the approach. We again use a similar case as in Section~\ref{sec:antithetic}, namely~\eqref{eq:1} with $\alpha = 1/2$ and $\beta = 3/2$. Considering only configurations that exactly satisfy~\eqref{eq:cond12_un}, we obtain the results shown on Figure~\ref{fig:SQS_un}. It is also possible, among the configurations that exactly satisfy~\eqref{eq:cond12_un}, to select configurations that satisfy as best as possible the condition~\eqref{eq:cond12_deux}. In practice, we generate 2000 configurations that exactly satisfy~\eqref{eq:cond12_un} and select among them the 100 configurations for which the difference between the left and the right-hand sides of~\eqref{eq:cond12_deux} is the smallest. We then obtain the results shown on Figure~\ref{fig:SQS_deux}.

In the case considered here (for which the contrast in the field $A$ is equal to 3), the variance is reduced by a factor 20 when using configurations that exactly satisfy~\eqref{eq:cond12_un}, and by a factor 300 if~\eqref{eq:cond12_deux} is enforced as well. To compare this variance reduction approach with the two previous ones, it is however needed to consider a case for which the contrast in $A$ is similar. In that case, the variance is reduced by a factor of 9 when using configurations that exactly satisfy~\eqref{eq:cond12_un}, and by a factor of 60 if~\eqref{eq:cond12_deux} is enforced as well.  

\medskip

In all the test cases we have considered (see~\cite{sqs,minvielle-these} for details), we have observed that the systematic error is kept approximately constant by the approach (it might even be reduced), while the variance is reduced by several orders of magnitude. Such an efficiency is achieved at almost no additional cost with respect to the classical Monte Carlo algorithm.

\section{Related issues and Further research}
\label{sec:future}

The studies we have reviewed above on different variance reduction approaches definitely show that such approaches may be very beneficial in the context of random homogenization, improving the accuracy while essentially preserving the computational cost. Their efficiency, measured as the actual ratio between the variance of a quantity computed with a direct Monte-Carlo approach and that of the same quantity computed using the variance reduction approaches, varies, depending upon the amount of information that one has on the problem and that one inserts into the specific variance reduction approach. The antithetic variable approach, a quite generic approach that can be put in action almost without  any prior knowledge on the problem considered, already reduces the variance by one order of magnitude, say, in the best case scenarios. Control variate and Special QuasiRandom Structures, both approaches that require exploiting some information on the problem, perform much better. Their efficiency may typically be one order of magnitude larger.

Of course, the efficiency of \emph{all} approaches is sensitive to the contrast present in the original multiscale problem. In a schematic manner, one may say that the efficiency is, approximately, inversely proportional to the contrast. It is an issue, since practically relevant multiscale problems may present a high contrast. Fortunately, there is room for improvement in the approaches and several ideas, some of them already explored in other contexts of the engineering sciences, some of them not, have not been pursued yet.

Among possible tracks for further research, we wish to cite a couple of alternate control variate approaches.

\medskip

A first possible track consists in considering {\em nonlinear} convex stochastic homogenization problems (as those considered in~\cite{dcds}), and use a corresponding {\em linear} problem either as a control variate (in the spirit of the approaches presented in Section~\ref{sec:controlvariate}) or as a way to select particular configurations (as in Section~\ref{sec:sqs}). We do not detail here the precise construction of this linear model, but rather focus on how to use it in practice. Let $\xi \in \RR^d \mapsto W^\star(\xi) \in \RR$ be the homogenized energy density of the nonlinear stochastic homogenization problem, and $\xi \mapsto W^\star_N(\omega,\xi)$ be its approximation computed by considering the nonlinear cell problem on the bounded domain $Q_N$. Let $A^\star_N(\omega)$ be the homogenized matrix of the corresponding linear problem. Our aim is to use $\xi^T A^\star_N(\omega) \xi$ as a control variate for $W^\star_N(\omega,\xi)$. Note however that we do not know the expectation of $\xi^T A^\star_N(\omega) \xi$, and hence we cannot directly use a Monte Carlo algorithm on the random variable
$$
W^\star_N(\omega,\xi) - \rho \Big( \xi^T A^\star_N(\omega) \xi - \EE \left[ \xi^T A^\star_N(\omega) \xi \right] \Big).
$$
However, computing $A^\star_N(\omega)$ is expected to be less expensive than computing $W^\star_N(\omega,\xi)$, because the corrector problem in the former case is linear, whereas it is nonlinear in the latter case. A natural idea is thus to replace, in the above relation, $\EE \left[ \xi^T A^\star_N(\omega) \xi \right]$ by an empirical mean. This leads to approximate $\EE \left[ W^\star_N(\omega,\xi) \right]$ by a mean of the form
$$
\frac{1}{M} \sum_{m=1}^M \Big( W^\star_N(\omega_m,\xi) - \rho \ \xi^T A^\star_N(\omega_m) \xi \Big) + \frac{\rho}{\cal M} \sum_{m=1}^{\cal M} \xi^T A^\star_N(\omega_m) \xi,
$$
where $M$, ${\cal M}$ (that we expect to be much larger than $M$) and $\rho$ are chosen to minimize the variance of the approximation for a given computational cost.

\medskip

A second track for further research is to use the so-called {\em bounds}, that are routinely employed in Mechanics, in order to build a control variate approach. Given the computational cost for obtaining approximations of $A^\star$, practitioners indeed sometimes choose to avoid computing the actual homogenized coefficients (by solving~\eqref{eq:AstarN}--\eqref{eq:correcteur-random-N}) and concentrate on \emph{bounds} (namely the Reuss, Voigt, Hashin-Shtrikman bounds, \dots) on the homogenized matrix~$A^\star$. 

For the sake of illustration, let us briefly review the derivation of the so-called Voigt bound. We assume that the random coefficient $A$ is a symmetric matrix. This assumption is critically used in what follows, and more generally in the derivation of many bounds. Under this assumption, the matrix $A^\star_N(\omega)$, defined by~\eqref{eq:AstarN}, satisfies, for any $p$,
$$
p^T A^\star_N(\omega) p = \inf \left\{ \frac{1}{|Q_N|} \int_{Q_N} (p+\nabla v)^T A(\cdot,\omega) (p+\nabla v), \quad v \in H^1_{\rm per}(Q_N) \right\}
$$
and hence, by choosing $v=0$ in the above problem, we obtain that
$$
A^\star_N(\omega) \leq \frac{1}{|Q_N|} \int_{Q_N} A(\cdot,\omega).
$$
The average of $A(\cdot,\omega)$ over $Q_N$ hence provides an upper bound on $A^\star_N(\omega)$, which is the so-called Voigt bound. 

In the specific case of two-phase composite materials (made of two phases denoted $\mathcal A$ and $\mathcal B$), where the random coefficient is given, with obvious notations, by 
$$
A(x,\omega) = \chi(x,\omega) \, A + (1-\chi(x,\omega)) \, B,
$$
where $\chi$ is the characteristic function of the phase $\mathcal A$, more elaborate bounds have been proposed, including the so-called Hashin-Shtrikman bounds. We refer e.g. to~\cite{singapour} for more details. The idea we are currently pursuing is to use these bounds not as an approximation for $A^\star_N(\omega)$, but as a control variate.

\enlargethispage{20pt}





\funding{The work of the last two authors is partially supported by EOARD under Grant FA8655-13-1-3061 and by ONR under Grant N00014-12-1-0383.}

\ack{The authors would like to thank all their collaborators on the issues presented here and related issues, in particular W.~Minvielle (Ecole des Ponts and INRIA).}



\bibliographystyle{plain}

%
%
%
%
%

\end{document}